\documentstyle[12pt]{article}

\begin{document}

\begin{center}

{\bf  Sensitive dependence and dense periodic points }\\
\vspace{.2in}
 S.Kanmani\\ 
 Materials  Science Division, Indira Gandhi Centre for Atomic Research, Kalpakkam- 603 102. India .

\end{center}

{ \bf Abstract- }  We show sensitive dependence on initial conditions and dense set of periodic points 
imply asymptotic sensitivity, a stronger form of sensitivity, where the deviation happens
not just once but infinitely many times. As a consequence it follows that all Devaney chaotic
systems ( e.g. logistic map ) have this asymptotic sensitivity.

\vspace{.2in}

{\bf Keywords-}  sensitivity, Devaney's chaos, dense periodic points, asymptotic sensitivity 
\vskip 0.1cm
	Sensitive dependence on initial conditions (shortly, sensitivity) is a
central concept in the theory chaos and  discrete dynamical systems.
Roughly speaking, a sensitive system with a sensitivity constant $ \delta >0 $ has
the following property. Arbitrarily close to the initial point of any specified
trajectory, there is a point whose trajectoty deviates from
that  specified one   by  a distance more than $\delta $, atleast at one instant of time later. 
The words  arbitrarily close  emphasizes the fact that even by  choosing the initial point
closer and closer to that specified inital point one cannot avoid a  deviation  in
  distance more than $  \delta $. Thus, the  arbitrarily small initial separation grows
up to  more than $ \delta > 0 $, within a finite time. Hence, in predicting the future
an error of magnitude $\delta $ is inevitable, however small be the error in the initial
condition.  This is often refered to, rather crudely, as the divergence of nearby trajectories.
 The following definition of sensitivity makes this precise.\\

	Discrete dynamical system is a pair  $( X,T )$ where X is a metric  
space and $T$ is a continuous self-map on $ X$. $B_{r}(x) $ denotes the open ball centered at $x$
and of radius $r$.  The distance between the points $x$ and $y$ of the space $X$ is denoted by $ d(x,y) $.
 The space $ X$ is supposed to represent all possible states of a physical system and $ T$ models the 
evolution of the system from the present  state   $x $ to the next state $ T(x) $. Thus the sequence 
$ \{ x, T(x), T^2(x), T^3(x),...  \} $ represents the evolution of the state $x$ and is called
the trajectory of $x$. If $T^k(x)=x $ for some natural number $ k$ then $x$ is periodic point.
The smallest $k$ for which  $T^k(x)=x $ is the order of the periodic point $x $.\\

{\bf Definition: (sensitivity) } The system  $(X,T)$ is said to be sensitive if there exists
a positive constant $ \delta $, called sensitivity constant, such that for all $ x \in X $ for all
$ \epsilon  > 0 $, there is a natural number $ n > 0 $ and a  pair of points $ y $ and $ z $ in
 $ B_{ \epsilon}(x)    $ such that $ d(T^n(z), T^n(y)) > \delta. $ \\

Remark: Note that if $ d(T^n(z),T^n(y)) > \delta $ then either $ d(T^n(x),T^n(z)) > \delta/2 $ or 
$ d(T^n(x),T^n(y)) > \delta/2 $. So equivalently with respect to the center of ball $B_{\epsilon}(x)$, there  is
a point ($ y$ or $z$ ) which deviates by $ \delta/2 $.  Some prefer this $ \delta/2$ version of the
definition.\\
The question that  we ask arises in  the following way. The above definition does not preclude 
the following possibility. All  those trajectories, that deviate to a distance more than $\delta $
 say at some  $k$-th instant, could get closer and closer to that specified trajectory ever after.
In other words, all these temporarliy deviant trajectories, could  asymptotically approach the
sepecifed trajectory. In such a system, eventually you will be making lesser
and lesser error, in  predicting the state of the  specified trajectory. Does this happen in
the standard sensitive systems like the logistic or a tent  map ?  Are there  sensitive maps
 in which the above scenario of asymptotic
convergence of  nearby trajetories is impossible ? In this note we answer these questions.
Any system which has the following property cannot have such asymptotic convergence of
nearby trajectories.\\
{\bf Definition: (Asymptotic sensitivity) }
The dynamical system $(X, T ) $ is said to be asymptotically sensitive if there is a $ \delta > 0 $
such that  any ball $ B_{ \epsilon}(x) $  ( arbitrary $ x \in X $ and arbitrary $\epsilon > 0 $) 
contains  points $ y $ and $ z $ such that  $ d(T^{n_{i}}(y), T^{n_{i}}(z)) > \delta $ for infinitely many
positive numbers $ n_{i} $.\\

Note: Consider the space $ X $ of one sided infinite strings from the alphabet $ \{ 0,1\} $. 
The  shift dynamical system S, acting on this space  maps  the string $ a = a_{0} a_{1} a_{2}....$  to    $S (a)=
a_{1}a_{2}a_{3}....  $  where $a_{i}$ is either 1 or 0. The distance between two strings $ a$ and $b$ is 
 $ d(a,b)= \sum_{i=0}^{\infty} \frac{D(a_{i},b_{i})}{2^{i}} $ where $ D$ is the discrete metric 
in the space of two symbols $ \{0,1 \}$. Let $ A$ be the subspace of $ X$ consisting of  strings whose 
 terms are eventually all zero i.e. beyond some $n$th term all the terms are zero.  The restriction of $T $ to $ A $ is an
example of a dynamical system which is sensitive but not asymptotically sensitive [1].\\
 {\bf Theorem:} If a dynamical system $(X,T)$ is sensitive and has dense set of periodic points then it is asymptotically sensitive.\\
 {\bf Proof:} Let the sensitivity constant of the system be $ \delta $. 
 The main idea is to choose the diverging pair of points to be periodic points. So one has to show that in any open ball,
there exists a pair of periodic point themselves, which separate to a distance more than $\delta$.  Since these points are
periodic, this separation would occur infinite number of timesleading to asymptotic sensitivity. By sensitivity, in any open
ball $B_{r}(x)$ there are points $z,y$ such that $$d(T^k(z), T^k(y)) >\delta $$ for some integer $ k > 1$.  Since $T^k$ is
continuous, for any open ball $ B_{\epsilon}(T^k(z))$ there is a $\rho >0 $ such that $$ T^k(B_{\rho }(z)) \subset
B_{\epsilon}(T^k(z))$$ since $ U= B_{r}(x) \cap B_{\rho}(z)$ is a non-empty open set (contains $z$ ), it has a periodic
point $ p$ (set of periodic points being dense ) such that $ T^k(p) \in B_{\epsilon}(T^k(z))$.  Similarly there is a
periodic point $q \in B_{r}(x)$ such that $ T^k(q) \in B_{\epsilon}(T^k(y)) $.  Now we have two periodic points $ p,q \in
B_{r}(x) $ such that $$ T^k(p) \in B_{\epsilon }(T^k(z)) $$ and $$ T^k(q) \in B_{\epsilon }(T^k(y)) $$ Since $\epsilon$ is
arbitrary it is possible to choose the $ \epsilon $ such that $ d(T^k(p),T^k(q)) > \delta $.  Let $\epsilon $ be chosen such
that $$ d(T^k(z),T^k(y))  -4\epsilon > \delta $$ Then applying triangle inequality to the points $ T^k(z), T^k(y)$ and $
T^k(p)$ one gets $$ d(T^k(p), T^k(y)) \geq d(T^k(z),T^k(y)) -d(T^k(z),T^k(p)) $$ and as $d(T^k(z),T^k(p)) < \epsilon $ it
follows $$ d(T^k(p),T^k(y)) \geq d(T^k(z),T^k(y)) - \epsilon $$ Similarly applying triangle inequality to $ T^k(p), T^k(q)$
and $ T^k(y)$ and using the fact that $d(T^k(q),T^k(y)) < \epsilon $ one gets $$d(T^k(p),T^k(q)) \geq d(T^k(p), T^k(y)) -
\epsilon $$ Combining these one gets $$d(T^k(p),T^k(q)) \geq d(T^k(z), T^k(y)) - 2\epsilon $$ As $ d(T^k(y), T^k(z)) >
\delta + 4\epsilon $ by the earlier choice of $\epsilon $ it follows $$ d(T^k(p), T^k(q)) > \delta $$ By
construction the points $p$ and $q$ are periodic points in the ball $B_{r}(x)$ and hence  $ d(T^{k_{i}}(p),
T^{k_{i}}(q)) >\delta $ for infinitely many positive integers $k_{i}$. This proves the proposition.\\

 Du [2],  proves amongst others that if the system $ ( X,T )$ is topologicaly transitive and has dense set of
periodic points then $ T^{k} $ is asymptotically sensitive for every positive $k$. However, we assume here only
sensitivity and dense set of periodic points. Example of a dynamical system which is sensitive
and has a dense set periodic points but is not topologically transitive is given in [3]. \\

Since Devaney chaotic systems are  sensitive and have dense set of periodic points it follows from the
theorem above that they are all asymptotically sensitive as well. For e.g. the well known dynamical
systems $( [0,1], 4x(1-x)) $ ( logisitic map) and $ ([0,1], 1-|2x-1| )$  (tent map ) are asymptotic sensitive.\\

{\bf Acknowledgement.} The author thanks Prof.V.Kannan and his open neighbourhood for many
illuminating discussions and Dr K.P.N. Murthy for his encouragement and support.

{\bf References}

1. V.Kannan, University of Hyderabad, Personal communication, (2002).

2. Bau sen Du, A dense orbit almost implies sensitivity to initial conditions, Bull. Inst. Math. Acad. Sinica.
26 (1998), 85-94. 

3. David Asaaf and Steve Gadbois, Definition of chaos, American mathematical Monthly, (1992) p-865.

 \end{document}